\documentclass[12pt]{article}

\newtheorem{theorem}{Theorem}[section]
\newtheorem{proposition}[theorem]{Proposition}
\newtheorem{lemma}[theorem]{Lemma}
\newtheorem{corollary}[theorem]{Corollary}

\newtheorem{example}[theorem]{Example}

\newcommand{\R}{{\bf R}}
\makeatletter
\@addtoreset{equation}{section}
\makeatother

\begin{document}

\title{Analytic Projective Finsler Metrics and Sprays}
\author{Zhongmin Shen}
\date{}
\maketitle

\begin{abstract} We study projective Finsler metrics and sprays on an open subset $U$ of $R^n$. They are characterized by 
the Hamel equations
\[
 P_{x^m y^k}y^m=P_{x^k},
\]
where $P=P(x, y)$ are function on $TU$ with positive homogeneity: $P(x, cy)= c (x,y)$, $c >0$. This is the Hilbert's Fourth problem in the regular case. We determine all solutions analytic in the base variable.   We compare this local representation with earlier work on Hilbert's fourth problem and with the restricted Taylor expansions arising for projective  Finsler metrics of constant flag curvature.   We conclude with explicit nonlinear projective Finsler  metrics and  sprays on Euclidean space, including one-sided complete examples and a nonflat complete example.
\end{abstract}

\section{Introduction}

Hilbert's fourth problem asks for the characterization of metrics on an open set of  $R^n$ whose length-minimizing curves are straight lines. Such metrics are called {\it projective metrics}. Regular metrics are called Finsler metrics. In his classical paper of 1903, Hamel studied projective Finsler metrics  and derived differential equations characterizing those metrics \cite{Hamel}. Hamel proved that 
for a Finsler metric $F=F(x,y)$ on an open set $U\subset\R^n$, its geodesics are straight lines  if and only if $F$ satisfies
\begin{equation}
 F_{x^m y^k}y^m=F_{x^k},\qquad k=1,\ldots,n.   \label{HamelEqF}
\end{equation}
In this case, the spray of $F$ is given by $G^i=P(x,y)y^i$, where $P=P(x,y)$ is called  the projective factor given by 
\[ P= \frac{F_{x^m}y^m}{2 F}.\]
The subject was later developed by methods of convex and integral geometry.  Pogorelov obtained a far-reaching solution of Hilbert's fourth problem and presented it systematically in \cite{Pogorelov}.
 Szabo  in his  paper  \cite{Szabo}, he develops a general framework for describing projective Finsler metrics, extending earlier work by Pogorelov.
Alvarez Paiva related projective Finsler metrics to symplectic geometry on the space of oriented lines \cite{AlvarezPaiva}, and Crampin explained the relation of this approach with the classical formulations of Hamel and Rapcs'ak \cite{Crampin}.

We may consider a more general structure - spray  $G$ on an open set of $R^n$. Its geodesics are straight lines if and only if $G^i= P(x,y)y^i$. 
Further $G$ is of isotropic curvature if and only if $P=P(x,y)$ satisfies the Hamel equations (\ref{HamelEqF})
\begin{equation}
P_{x^m y^k}y^m=P_{x^k},\qquad k=1,\ldots,n.   \label{HamelEqP}
\end{equation}
The solutions of (\ref{HamelEqP}) are called {\it Hamel functions}. 

A natural question is how to study the Hamel functions.
The aim of the present paper is to give a self-contained coefficientwise answer in the $x$-analytic category.   Inspired by the author's earlier work on projective flat Finsler metrics of constant flag curvature \cite{Shen}, we   shall show that if a Hamel function  $P=P(x,y)$ is analytic in $x$, then it can expressed by
\begin{equation}
 P=\sum_{m=0}^{\infty}{1\over m!}
 (x^i\partial_{y^i})^m\Phi_m(y) =\sum_{m=0}^{\infty} \frac{1}{m!} \frac{d^m}{dt^m} \Big [ \Phi_m( y+tx) \Big ]|_{t=0},\label{eq1.3}
\end{equation}
where each $\Phi_m$ is smooth on the slit fiber and positively homogeneous of degree $m+1$. Conversely, any series of the form (\ref{eq1.3}) that converges together with the derivatives required in the Hamel equations defines a solution. See Theorem \ref{thm4.1} below.

This result shows that the general analytic Hamel solution is parametrized by the sequence
\[
 \Phi_0,\Phi_1,\Phi_2,\ldots,
\]
where $\Phi_m$ has degree $m+1$. These functions are called the {\it potentials} of $P$.

If $F=F(x,y)$ is a Finsler metric satisfying (\ref{HamelEqF}), then 
 we obtain a Taylor expansion for  $F$:  
\begin{equation}
F=\sum_{m=0}^{\infty}{1\over m!}
 (x^i\partial_{y^i})^m\Psi_m(y) =\sum_{m=0}^{\infty} \frac{1}{m!} \frac{d^m}{dt^m} \Big [ \Psi_m( y+tx) \Big ]|_{t=0},\label{eq1.4}
\end{equation}
where each $\Psi_m$ is smooth on the slit fiber and positively homogeneous of degree $m+1$. However, the projective factor of $F$ is not necessarily a Hamel function unless $F$ is of isotropic curvature. In this case, the potentials of $F$ and its projective factor $P$ take special forms.  See Example \ref{ex4.1} below.

We shall construct several interesting projective sprays of isotropic curvature. For example, if taking 
$\Phi_2 (y) = \frac{2}{3} |y|^3$ and $\Phi_m (y) =0 $ for $m \not=2$, we get 
\[ P(x, y) = |x|^2 |y| +\frac{\langle x, y \rangle^2}{ |y|}.\]
This simple projective spray is forward complete with isotropic curvature. 
Note that $P=P(x,y) \geq 0 $ is not a Finsler metric. Adding  it to the standard Euclidean metric, we obtain a projective Finsler metric  $ F:= |y| + P (x,y)$  on $R^n$. However, the projective factor  of $F$ is no longer a Hamel function since $F$ is not of isotropic flag curvature.

\section{Finsler metrics, sprays, and isotropic curvature}

Let $M$ be a smooth $n$-dimensional manifold, and let $U\subset M$
be a coordinate neighborhood with local coordinates $(x^i)$.
The induced standard local coordinates on $TU$ are $(x^i,y^i)$,
where $y=y^i\partial/\partial x^i$.  
A Finsler metric on $M$ is a
smooth function $F:TM\setminus\{0\}\to(0,\infty)$ that is positively
homogeneous of degree one in $y$ and whose fundamental tensor
\begin{equation}\label{eq:fundamental-tensor}
 g_{ij}(x,y)={1\over2}(F^2)_{y^iy^j}(x,y)
\end{equation}
is positive definite.  The geodesics of $F$ are characterized by the following equations:
\begin{equation}\label{eq:finsler-geodesics}
 \ddot x^i+2G^i(x,\dot x)=0,
\end{equation}
where
\begin{equation}\label{eq:finsler-spray}
 G^i={1\over4}g^{ij}\left\{
 (F^2)_{x^ky^j}y^k-(F^2)_{x^j}\right\}.
\end{equation}
The local functions $G^i$ are positively homogeneous of degree two in $y$. These local functions $G^i$ define a vector field $G$ on $TM$ 
whose expression in the standard local coordinates $(x^i,y^i)$ is
\begin{equation}\label{eq:general-spray-vector-field}
 { G}=y^i{\partial\over\partial x^i}
       -2G^i(x,y){\partial\over\partial y^i}.
\end{equation}

More general, one can consider a vector field $G$ locally expressed in the form (\ref{eq:general-spray-vector-field}).
 Its integral curves project to curves on the manifold
satisfying the geodesic equations (\ref{eq:finsler-geodesics}). This kind of vector fields are called {\it sprays}. 
Sprays are not necessarily induced by any (pseudo-) Finsler metric. 

The Riemann curvature (Jacobi endomorphism) of a spray is the family
of functions
\begin{equation}\label{eq:spray-riemann-curvature}
 R^i{}_k=
 2{\partial G^i\over\partial x^k}
 -y^j{\partial^2G^i\over\partial x^j\partial y^k}
 +2G^j{\partial^2G^i\over\partial y^j\partial y^k}
 -{\partial G^i\over\partial y^j}
  {\partial G^j\over\partial y^k}.
\end{equation}
Its Ricci curvature is $\mathrm{Ric}=R^i{}_i$.  For $n\geq2$ put
$R=\mathrm{Ric}/(n-1)$.  We say that a spray has
{\it isotropic curvature} if
\begin{equation}\label{eq:isotropic-curvature}
 R^i{}_k=R\delta^i_k-{1\over2}R_{y^k}y^i.
\end{equation}
Here $R$ is positively homogeneous of degree two in $y$.

A Finsler metric is called {\it projective} (or projectively flat)
in the chosen coordinate neighborhood if its geodesics, regarded as
unparameterized curves, are straight lines in the coordinates $(x^i)$.  By Hamel's characterization, this
holds precisely when  $F$ satisfies the Hamel equations (\ref{HamelEqF}): 
\begin{equation}\label{eq:finsler-hamel-section}
 F_{x^ky^j}y^k=F_{x^j}.
\end{equation}
In this case the induced spray (\ref{eq:finsler-spray}) has the form
\begin{equation}\label{eq:finsler-projective-factor}
 G^i=P(x,y)y^i,\qquad
 P(x,y)={F_{x^k}(x,y)y^k\over2F(x,y)}.
\end{equation}
Conversely, if the spray induced by a Finsler metric has the form
$G^i=P(x,y)y^i$, then its unparameterized geodesics are straight lines,
so that $F$ satisfies (\ref{eq:finsler-hamel-section}).

A spray $G=y^i{\partial\over\partial x^i}
       -2G^i(x,y){\partial\over\partial y^i}$ on $M$ is called {\it projective } in these coordinates if
\begin{equation}\label{eq:projective-spray-definition}
 G^i(x,y)=P(x,y)y^i.
\end{equation}
For convenience write $P_0=P_{x^j}y^j$.  Substituting
(\ref{eq:projective-spray-definition}) into
(\ref{eq:spray-riemann-curvature}) gives
\begin{equation}\label{eq:projective-spray-curvature}
 R^i{}_k=(P^2-P_0)\delta^i_k
 +\big(3P_{x^k}-(P_0)_{y^k}-PP_{y^k}\big)y^i.
\end{equation}
In particular,
\begin{equation}\label{eq:projective-spray-scalar}
 R=P^2-P_0,\qquad \mathrm{Ric}=(n-1)(P^2-P_0).
\end{equation}
The characterization of Li and Shen \cite{LiShen} takes the
following form.

\begin{proposition}\label{prop:li-shen-hamel} (\cite{LiShen})
A projective spray $G^i=P(x, y)y^i$ on $U$ has isotropic curvature if
and only if its projective factor $P$ satisfies the Hamel equations
\begin{equation}\label{eq:projective-hamel-section}
 P_{x^jy^k}y^j=P_{x^k}.
\end{equation}
\end{proposition}

\noindent{\bf Proof.}
Since
\[
 (P_0)_{y^k}=P_{x^k}+P_{x^jy^k}y^j,
 \qquad R_{y^k}=2PP_{y^k}-(P_0)_{y^k},
\]
the coefficient of $y^i$ in
$R^i{}_k-R\delta^i_k+\frac12R_{y^k}y^i$ is
\[
 3P_{x^k}-(P_0)_{y^k}-PP_{y^k}
 +\frac12\big(2PP_{y^k}-(P_0)_{y^k}\big)
 ={3\over2}\big(P_{x^k}-P_{x^jy^k}y^j\big).
\]
Thus (\ref{eq:isotropic-curvature}) is equivalent to
(\ref{eq:projective-hamel-section}). \hfill{\bf Q.E.D.}

Note that (\ref{eq:finsler-hamel-section}) and (\ref{eq:projective-hamel-section}) both are Hamel functions.

\begin{corollary}\label{cor:projective-finsler-isotropic}
Let $F$ be a Finsler metric on an open subset $U$ and put
$P=F_{x^j}y^j/(2F)$.  Then $F$ is projective and its induced
spray has isotropic curvature if and only if both $F$ and $P$
satisfy the Hamel equations:
\[
 F_{x^jy^k}y^j=F_{x^k},\qquad
 P_{x^jy^k}y^j=P_{x^k}.
\]
\end{corollary}

\bigskip

\bigskip

\noindent{\bf Proof.}
The first equation characterizes projective of $F$ and
identifies $P$ as its projective factor by
(\ref{eq:finsler-projective-factor}).  The second equation is
equivalent to isotropic curvature of the induced projective spray
by Proposition \ref{prop:li-shen-hamel}. \hfill{\bf Q.E.D.}

\section{Analytic solutions}

In this section, we shall give an analytic solution to the Hamel equations (\ref{eq:projective-hamel-section}).                

We first consider a single homogeneous polynomial in $x$.

\begin{lemma} \label{lem3.1}
Let
\begin{equation}
 P_m(x,y)={1\over m!}A_{i_1\cdots i_m}(y)
 x^{i_1}\cdots x^{i_m},
\end{equation}
where $A_{i_1\cdots i_m}$ is symmetric in its indices and homogeneous of degree one in $y$.  Then $P=P_m(x, y)$ satisfies the Hamel equations 
(\ref{eq:projective-hamel-section}) if and only if there is a function $\Phi_m(y)$, homogeneous of degree $m+1$, such that
\[
 A_{i_1\cdots i_m}
 =(\Phi_m)_{y^{i_1} y^{i_2} y^{i_2}\cdots y^{i_1}},
\]
that is, $P_m$ is given by
\begin{equation}
P_m(x,y) = \frac{1}{m!} (x^i\partial_{y^i})^m \Phi_m (y).
\end{equation}
\end{lemma}

\noindent{\bf Proof.}
For $m=0$ the assertion is immediate.  Assume $m\geq1$.  Substitution of the defining expression for $P_m$ into the Hamel equations
\[ ({P_m})_{x^j y^k} y^j = (P_m)_{x^k}\]
 gives
\begin{equation}
 y^j(A_{j i_2\cdots i_m})_{y^k}
 =A_{k i_2\cdots i_m}.
\end{equation}
Define
\begin{equation}
 \Phi_m (y) = \frac{1}{m+1} P_m(y, y) = \frac{1}{(m+1)!} A_{j_1\cdots j_m}(y)y^{j_1}\cdots y^{j_m}.
\end{equation}
We claim that for $1\leq r\leq m$,
\begin{equation}
(\Phi_m)_{y^{i_1}\cdots\partial y^{i_r}}
 ={1\over(m+1-r)!}
 A_{i_1\cdots i_rj_{r+1}\cdots j_m}
 y^{j_{r+1}}\cdots y^{j_m}.
\end{equation}

For $r=1$, differentiating the definition of $\Phi_m$ gives one term in which $A$ is differentiated and $m$ terms in which one of the explicit factors $y^{j_a}$ is differentiated.  By symmetry and the coefficient identity above,
\[
 (\Phi_m)_{y^k}=\frac{1}{m!} A_{k j_2\cdots j_m}
 y^{j_2}\cdots y^{j_m},
\]
which is the claimed formula for $r=1$.

Suppose the claimed formula holds for some $r<m$.  Differentiate once more.  Differentiating the $m-r$ explicit $y$ factors gives $m-r$ identical terms.  Differentiating $A$ and contracting the derivative with one of the remaining $y$ factors gives, by the coefficient identity and the symmetry of $A$, one additional identical term.  Thus the new factor is
\[
 (m-r)+1=m+1-r.
\]
This proves the formula by induction.

Putting $r=m$ gives
\begin{equation}
 (\Phi_m)_{y^{i_1}\cdots y^{i_m}}
 =A_{i_1\cdots i_m}.
\end{equation}
Consequently, $\Phi_m$ has the required derivatives.  Since $A$ is one-homogeneous, $\Phi_m$ are homogeneous of degree $m+1$.

Conversely, suppose $\Phi_m$ is $(m+1)$-homogeneous and define $A$ by the stated derivative formula.  Then
\[
 (\Phi_m)_{y^k y^{i_2}\cdots y^{i_m}}
\]
is homogeneous of degree one.  Euler's identity yields
\[
 y^j(\Phi_m)_{y^j y^k y^{i_2}\cdots y^{i_m}}
 =(\Phi_m)_{y^k y^{i_2}\cdots y^{i_m}},
\]
which is precisely the coefficient form of the Hamel equations.  Thus $P_m$ satisfies the Hamel equation. \hfill{\bf Q.E.D.}

We now obtain the analytic version.

\begin{theorem}\label{thm4.1}
Let $P(x,y)$ be smooth on a neighborhood of $\{0\}\times(\R^n\setminus\{0\})$, positively homogeneous of degree one in $y$, and real analytic in $x$ near $x=0$, locally uniformly with the mixed derivatives used below.  Then $P$ satisfies the Hamel equations if and only if it admits a locally convergent expansion
\begin{equation}
 P(x,y)=\sum_{m=0}^{\infty}{1\over m!}
 (x^i\partial_{y^i})^m\Phi_m(y), \label{eq4.1}
\end{equation}
where the series and its derivatives through one $x$-derivative and one $y$-derivative converge locally uniformly, and
\[
 \Phi_m(cy)=c^{m+1}\Phi_m(y),\qquad c>0.
\]
For a given analytic solution, the functions $\Phi_m$ are uniquely determined by its Taylor coefficients.
\end{theorem}

\noindent{\bf Proof.}
Expand $P$ at $x=0$:
\begin{equation}
 P(x,y)=\sum_{m=0}^{\infty}{1\over m!}
 A_{i_1\cdots i_m}(y)x^{i_1}\cdots x^{i_m},
\end{equation}
where
\[
 A_{i_1\cdots i_m}(y)
 =P_{x^{i_1}\cdots x^{i_m}}(0,y).
\]
The tensors $A$ are symmetric and one-homogeneous in $y$.  Comparing coefficients in the Taylor expansion of the Hamel equation shows that each homogeneous Taylor term satisfies the same equation. This comparison is justified by the stated analyticity and regularity assumptions.  Lemma \ref{lem3.1} therefore supplies a unique $(m+1)$-homogeneous $\Phi_m$ satisfying
\[
 A_{i_1\cdots i_m}=(\Phi_m)_{y^{i_1}\cdots y^{i_m}}.
\]
This gives (\ref{eq4.1}).

Conversely, every individual term in (\ref{eq4.1}) satisfies the Hamel equations by Lemma \ref{lem3.1}.  Whenever the series and the derivatives involved converge locally so that termwise differentiation is valid, the sum satisfies the Hamel equations. \hfill{\bf Q.E.D.}

\bigskip

Let $P=P(x,y)$ be a Hamel function on $U\subset R^n$ given by (\ref{eq4.1}). 
Then
\begin{equation}
 P(ty,y)=\sum_{m=0}^\infty (m+1)t^m\Phi_m(y).\label{PPhi_m}
\end{equation}
We obtain
\begin{equation}
 \Phi_m(y)={1\over (m+1)!}
 \left.{d^m\over dt^m}P(ty,y)\right|_{t=0}.
\end{equation}
This formula is useful to find the potentials $\Phi_m$ for a Hamel function $P$.

\section{Examples}

For consistency with Theorem \ref{thm4.1}, we specify the homogeneous
potentials for each example.

\begin{example}\label{ex4.1}({\bf Projective Finsler Metrics of Constant Flag Curvature} (\cite{Shen}))
Let $F=F(x, y)$ be a projective Finsler metric of constant flag  curvature ${\bf K} 
= \lambda $ on an open neighborhood  ${\cal U} $ of the origin in $ \R^n$. Suppose that 
$F$ is $x$-analytic at $x=0$. Then both
$F$ and  and its projective factor $P$  are given by 
\begin{eqnarray}
F(x, y ) & = &  \sum_{m=0}^{\infty} 
{1\over m!} {d^m \over dt^m }
\Big [ \Psi_m (y+tx) \Big ]\Big |_{t=0}\label{thm2F}\\
P(x, y) & = &  \sum_{m=0}^{\infty} 
{1\over m!} {d^m \over dt^m }
\Big [ \Phi_m (y+tx) \Big ]\Big |_{t=0},  \label{thm2P}
\end{eqnarray}
where
\[\Psi_m := \cases{ {1\over 2(m+1)\sqrt{-\lambda}} \Big \{
\Big (  \varphi+\sqrt{-\lambda}\psi \Big )^{m+1}- \Big (  \varphi-\sqrt{-\lambda}\psi \Big )^{m+1} \Big \} \ & if $\lambda \not=0 $ \cr\\
& \cr \\
 \psi \varphi^m \; & if $\lambda =0$ \cr}
\]
\smallskip
\[
\Phi_m :=\cases{ {1\over 2(m+1)} \Big \{
\Big (  \varphi+\sqrt{-\lambda}\psi \Big )^{m+1}+ \Big (  \varphi-\sqrt{-\lambda}\psi \Big )^{m+1} \Big \} \ & if $\lambda \not=0 $ \cr\\
& \cr \\
 \frac{1}{m+1}\varphi^{m+1} \; & if $\lambda =0$ \cr}
\]
where $\psi(y):=F(0, y)$ and $\varphi(y):=P(0, y)$.
Conversely,
for any Minkowski norm $\psi(y)$ on $\R^n$, and any positively homogeneous function of degree one $\varphi(y)$ on $\R^n$,  the function $F(x,y)$ defined in (\ref{thm2F}) is a projective Finsler metric with constant flag  curvature ${\bf K}=\lambda$ and its projective factor $P(x,y)$ is given by (\ref{thm2P})
provided that they are convergent.

\medskip
\noindent{\bf The Riemannian case (Beltrami--Klein model).}
Consider the Riemannian metric of constant sectional curvature
$\lambda$ in projective coordinates,
\begin{equation}\label{eq:klein-metric}
 F_\lambda(x,y)=
 {\sqrt{(1+\lambda|x|^2)|y|^2
              -\lambda\langle x,y\rangle^2}
  \over 1+\lambda|x|^2}.
\end{equation}
For $\lambda<0$ we take $|x|<1/\sqrt{-\lambda}$; for
$\lambda\geq0$ the formula is defined on $\R^n$.
Its Riemannian tensor is
\[
 g_{ij}(x)={ (1+\lambda|x|^2)\delta_{ij}
                 -\lambda x_i x_j
              \over (1+\lambda|x|^2)^2}.
\]
The metric has constant sectional curvature $\lambda$ and its
unparameterized geodesics are straight lines.  More precisely,
its Christoffel symbols are
\[
 \Gamma^i_{jk}=-{\lambda\over1+\lambda|x|^2}
                    (x_j\delta^i_k+x_k\delta^i_j),
\]
so that its spray coefficients are $G^i=P_\lambda y^i$, where
\begin{equation}\label{eq:klein-projective-factor}
 P_\lambda(x,y)=-{\lambda\langle x,y\rangle
                         \over1+\lambda|x|^2}.
\end{equation}
In particular, $P_\lambda$ is linear in $y$, whereas
$F_\lambda$ is a nonlinear Hamel function.

We now determine the potentials in the Taylor representation
(\ref{thm2F})--(\ref{thm2P}).  Along $x=ty$ we have
\[
 F_\lambda(ty,y)={|y|\over1+\lambda t^2|y|^2},
 \qquad
 P_\lambda(ty,y)=-{\lambda t|y|^2
                         \over1+\lambda t^2|y|^2}.
\]
Comparing these geometric series with (\ref{PPhi_m}) gives,
for $k=0,1,2,\ldots$,
\begin{equation}\label{eq:klein-F-potentials}
 \Psi_{2k}(y)={(-\lambda)^k\over2k+1}|y|^{2k+1},
 \qquad \Psi_{2k+1}(y)=0,
\end{equation}
\begin{equation}\label{eq:klein-P-potentials}
 \Phi_{2k}(y)=0,
 \qquad
 \Phi_{2k+1}(y)={(-\lambda)^{k+1}\over2k+2}|y|^{2k+2}.
\end{equation}
These formulas also follow directly from the general expressions
above by taking $\psi(y)=|y|$ and $\varphi(y)=0$.
For $\lambda=0$ the metric is Euclidean, $F_0=|y|$ and
$P_0=0$.  For $\lambda=-1$, (\ref{eq:klein-metric}) is the
usual Klein metric on the unit ball; for $\lambda>0$ it is the
gnomonic (projective) coordinate expression of the round metric
on a hemisphere.

\end{example}

\begin{example}[Funk-type deformations of the Klein metric]\label{ex:funk-type}
Let $\lambda\in\R$ and $c\in\R$ be constants.  Set
\[
 Q=\sqrt{(1+\lambda |x|^2)|y|^2-\lambda \langle x,y\rangle^2}.
\]
Starting from the Klein metric in (\ref{eq:klein-metric}), consider
\begin{equation}\label{eq:funk-type-metric}
 F_{\lambda,c}(x,y)=F_\lambda(x,y)+{c\langle x,y\rangle\over1+\lambda|x|^2}.
\end{equation}
The second summand is an exact one-form: for $\lambda\ne0$,
\[
 {c\langle x,y\rangle\over 1+\lambda |x|^2}
 ={c\over2\lambda}\,d_x\Big [ \log (1+\lambda |x|^2) \Big ] (y),
\]
whereas for $\lambda=0$ it is $d_x(c|x|^2/2)(y)$.
Thus $F_{\lambda,c}$ is projectively flat wherever it is a Finsler
metric.  In particular, it is a strongly convex Finsler metric near
$x=0$ for every real $\lambda,c$.  More explicitly, the Randers
condition $|\beta|_\alpha<1$ here is $|c|\,|x|<1$; we also require
$1+\lambda |x|^2>0$.  Consequently $\lambda>0$ is allowed on such a neighborhood,
even though the usual Funk metric is obtained for negative $\lambda$.

Since $F_\lambda$ has projective factor
$P_{\lambda,c}=-\lambda \langle x,y\rangle/(1+\lambda |x|^2)$, differentiation along $y$ yields
\begin{equation}\label{eq:funk-type-factor}
 P_{\lambda,c}(x,y)
 ={F_{\lambda,c,x^k}y^k\over2F_{\lambda,c}}
 ={cQ^2-2\lambda \langle x,y\rangle Q-c\lambda \langle x,y\rangle^2\over
   2(1+\lambda |x|^2)(Q+c\langle x,y\rangle)}.
\end{equation}
If $\lambda=-c^2$ and $c>0$, then
\begin{equation}\label{eq:funk-special}
 F_{-c^2,c}(x,y)=
 {\sqrt{(1-c^2|x|^2)|y|^2+c^2\langle x,y\rangle^2}+c\langle x,y\rangle
  \over1-c^2|x|^2},\qquad
 P_{-c^2,c}={c\over2}F_{-c^2,c}.
\end{equation}
This is the forward Funk metric on the ball $|x|<1/c$,
with constant flag curvature $-c^2/4$.
The choice $c<0$ gives the reverse Funk metric.

We determine all the homogeneous potentials in the expansions
(\ref{thm2F})--(\ref{thm2P}).  
On the line $x=ty$ one has $Q(ty, y)=|y|$, and hence
\begin{equation}\label{eq:funk-line-F}
 F_{\lambda,c}(ty,y)
 ={|y|+ct|y|^2\over1+\lambda t^2|y|^2}.
\end{equation}
The projective-factor formula also gives
\begin{equation}\label{eq:funk-line-P}
 P_{\lambda,c}(ty,y)
 ={1\over2}{d\over dt}\log F_{\lambda,c}(ty,y)
 ={c|y|\over2(1+ct|y|)}
  -{\lambda t|y|^2\over 1+\lambda t^2|y|^2}.
\end{equation}
By comparison with (\ref{PPhi_m}), the potentials of the metric
(denoted by $\Psi_m$) are, for $k=0,1,2,\ldots$,
\begin{equation}\label{eq:funk-F-potentials}
 \Psi_{2k}(y)={(-\lambda)^k\over2k+1}|y|^{2k+1},\qquad
 \Psi_{2k+1}(y)={c(-\lambda)^k\over2k+2}|y|^{2k+2}.
\end{equation}
In general, the projective factor $P_{\lambda,c}$ is not a Hamel function since $F_{\lambda,c}$ is not of isotropic flag curvature. 

For $\lambda=-c^2$, $ F_{-c^2,c}$ is of isotropic curvature, then $P_{-c^2, c}$ is a Hamel function with potentials 
\begin{equation}\label{eq:funk-special-potentials}
 \Phi_m(y)={c\over2}\Psi_m(y),\qquad m\geq0,
\end{equation}
in agreement with (\ref{eq:funk-special}).

\end{example}

\begin{example} [Exact one-form terms]\label{ex5.2}
Let $f=f(x)$ be smooth.  Then
\begin{equation}
 P(x,y)=f_{x^i}(x)y^i
\end{equation}
satisfies the Hamel equation, since
\[
 P_{x^m y^k}y^m=f_{x^m x^k}y^m=P_{x^k}.
\]
Thus
\[
 P=P_0(y)+df_x(y)
\]
is a basic family of solutions.  If $f$ is analytic at the origin, its potentials are
\[
 \Phi_m(y)={1\over(m+1)!}
 \left.{d^{m+1}\over dt^{m+1}}f(ty)\right|_{t=0}.
\]
The additional term $P_0(y)$ contributes only to $\Phi_0$.
The analytic theorem contains this family but also permits solutions
which are genuinely nonlinear in $y$.
\end{example}

\begin{example}[Degree one in $x$]
Let $H(y)$ be any two-homogeneous function.  Then
\begin{equation}
 P(x,y)=x^iH_{y^i}(y)
\end{equation}
satisfies Hamel equation.  
The only nonzero potential is $\Phi_1=H$; all $\Phi_m$ with
$m\ne1$ vanish.  Choosing a nonquadratic two-homogeneous $H$
gives a solution nonlinear in $y$.
\end{example}

\begin{example}[A quadratic solution in $x$]
Take
\[
 \Phi_2(y)={2\over3}|y|^3,\qquad \Phi_m=0\quad(m\ne2),
\]
where $|y|=(\sum_i(y^i)^2)^{1/2}$.  Then
\[
 (\Phi_2)_{y^iy^j}
 =2|y|\delta_{ij}+2{y_i y_j\over |y|}.
\]
Hence
\begin{equation}
 P(x,y)=|x|^2|y|+{\langle x,y\rangle^2\over |y|}
\end{equation}
satisfies the Hamel equations.  This is a simple explicit solution which is quadratic in $x$ and nonlinear in $y$.

\end{example}

\begin{example}[Elementary solutions generated by one linear form]
Let $a=(a_i)$ be a constant covector and let $f$ be an arbitrary smooth function of one variable.  Then
\begin{equation}
 P(x,y)=f(\langle a,x\rangle)\langle a,y\rangle
\end{equation}
satisfies the Hamel equation.  In particular, choosing
\[
 f(s)=e^s,\quad \sin s,\quad {1\over1-s},\quad \hbox{or}\quad \log(1+s)
\]
gives elementary closed-form solutions on their natural domains.  If $f$ is analytic at zero, then
\[
 \Phi_m(y)={f^{(m)}(0)\over(m+1)!}
                 \langle a,y\rangle^{m+1},\qquad m\geq0.
\]
For instance, $f(s)=e^s$ gives
\[ \Phi_m(y)=\frac{1}{(m+1)!}\langle a,y\rangle^{m+1}.\]
Finite sums with different constant covectors are again Hamel functions,
with potentials obtained by addition.
\end{example}

\begin{example}[Integral-geometric solutions]
Let $\mu$ be a measure on a suitable set of covectors $\xi$ and let $a(s,\xi)$ be smooth.  Formally,
\begin{equation}
 P(x,y)=\int a(\xi\cdot x,\xi)|\xi\cdot y|\,d\mu(\xi)
\end{equation}
satisfies the Hamel equations wherever differentiation under the integral is justified.  Indeed,
\[
 P_{x^k}=\int a_s(\xi\cdot x,\xi)\xi_k|\xi\cdot y|\,d\mu,
\]
while
\[
 P_{x^m y^k}y^m
 =\int a_s(\xi\cdot x,\xi)\xi_k (\xi\cdot y)
 {\rm sign}(\xi\cdot y)\,d\mu
 =P_{x^k}.
\]
When $a(s,\xi)$ is analytic in $s$ and the integral permits
termwise differentiation, its potentials are
\[
 \Phi_m(y)={1\over(m+1)!}
 \int \partial_s^m a(0,\xi)
       (\xi\cdot y)^m|\xi\cdot y|\,d\mu(\xi).
\]
These expressions are $(m+1)$-homogeneous.  Smoothness in $y$
requires suitable assumptions on the measure and integrand,
since the absolute value is not smooth on its zero set.
This family is related to the integral-geometric viewpoint in
Hilbert's fourth problem.
\end{example}

\section{Complete Projective Sprays}

The homogeneous-potential representation gives a simple way to construct projective
sprays with controlled global behavior.  Here the base domain is all of $\R^n$; whenever an infinite series is used, we assume the local convergence of all derivatives needed to define a smooth Hamel function.

\begin{proposition}[A sufficient condition on the potentials]\label{prop8.1}
Let $\Phi_m$ be positively homogeneous of degree $m+1$, and suppose
\[
 P(x,y)=\sum_{m=0}^{\infty}{1\over m!}
       (x^i\partial_{y^i})^m\Phi_m(y)
\]
is a smooth Hamel function on $\R^n\times(\R^n\setminus\{0\})$.
If  $P\geq 0$, 
then the projective spray $G^i=P(x,y)y^i$ is forward complete.
\end{proposition}

\noindent{\it Proof.}
For an oriented Euclidean line $x=x_0+sv$ with $v\ne0$, write
$p(s)=P(x_0+sv,v)$.  A geodesic of $G^i=Py^i$ is a reparametrization
of this line. Let $x^i (t) = x^i_0 + s(t) v^i$ be a geodesic. Plugging it into (\ref{eq:finsler-geodesics}) yields
\[\frac{d^2 s}{dt^2}  + 2 p(s) \big ( \frac{ds}{dt} \big )^2 =0.\]
Then
\[
\frac{ds}{dt}=C\exp\left(-2\int_0^s p(u)\,du\right),\qquad C>0.
\]
  By assumption, $p(s)\ge0$ for $s\ge0$, so $dt/ds\ge C$.
Thus $\lim_{s\to\infty} t =+\infty$. The spray is forward complete.

\hfill {\it Q.E.D.}

\begin{example}[A nonlinear entire Hamel function]
Choose a single nonzero homogeneous potential,
\[
 \Phi_2(y)={2\over3}|y|^3,\qquad \Phi_m=0\quad(m\ne2).
\]
The representation formula gives
\begin{equation}
 P(x,y)={1\over2}(x^i\partial_{y^i})^2\Phi_2(y)
       =|x|^2|y|+{\langle x,y\rangle^2\over|y|}.
\end{equation}
Thus $P$ is smooth for $y\ne0$, positively one-homogeneous,
quadratic in $x$, and genuinely nonlinear in $y$.  It satisfies the
Hamel equation by Theorem \ref{thm4.1}.  Since $P\ge0$, Proposition \ref{prop8.1}
shows that $G^i=P(x, y)y^i$ is forward complete on $\R^n$.

The example is not backward complete.  Along the line $x=sv$ with
$|v|=1$, one has $p(s)=2s^2$ and hence
\[
 {dt\over ds}=C\exp\left({4\over3}s^3\right).
\]
The integral of this expression over $(-\infty,0]$ is finite.
Consequently the same geodesic reaches Euclidean infinity in finite
negative affine time.

\end{example}

\begin{example}[A bounded exact potential]\label{ex8.2}
Let $f(x)=\exp(-|x|^2)$ and put
\[
 P(x,y)=df_x(y)=-2\exp(-|x|^2)\langle x,y\rangle.
\]
This is an entire analytic Hamel function, and its Taylor potentials
are explicitly
\[
 \Phi_{2k}(y)=0,\qquad
 \Phi_{2k+1}(y)={(-1)^{k+1}\over(k+1)!}|y|^{2k+2}
 \quad(k\geq0).
\]
Indeed, these follow from the Taylor expansion of
$f(ty)=\exp(-t^2|y|^2)$ and the recovery formula in Example \ref{ex5.2}.
Along any oriented line,
\[
 \int_0^sP(x_0+uv,v)\,du
       =f(x_0+sv)-f(x_0),
\]
which is bounded above and below.  Thus $dt/ds$ is bounded away
from zero in both directions and the associated spray  $G^i = P(x,y) y^i$ is complete
in both affine-time directions. 
\end{example}

\begin{example}[A sign-changing, one-sided complete Hamel function]
Fix $c>0$ and set,
\begin{equation}
 P(x,y)=u+c\left(|x|^2|y|+{\langle x,y\rangle^2\over |y|}\right).
\end{equation}
The nonzero potentials are $\Phi_1(y)=|y|^2/2$ and
$\Phi_2(y)=2c|y|^3/3$.  Thus $P$ is a globally defined Hamel
function, and it is genuinely nonlinear in $y$.  It changes sign:
for a unit vector $v$, $P(sv,v)=s+2cs^2$ is negative for
$-1/(2c)<s<0$ and positive for $s>0$.
At $x=0$, $P(0,y)=0$ and $P_{x^k}(0,y)y^k=|y|^2$; hence
\[
 P^2-P_{x^k}y^k=-|y|^2\ne0
 \qquad\hbox{at }x=0.
\]
For any oriented line $x=a+sv$, $v\ne0$, the leading term of
$p(s)=P(a+sv,v)$ is $2c|v|^3s^2$.  Consequently
\[
 2\int_0^s p(\sigma)\,d\sigma
       ={4c|v|^3\over3}s^3+O(s^2).
\]
The affine-time density $dt/ds$ grows super-exponentially as
$s\to+\infty$ and is integrable as $s\to-\infty$.  Thus the
associated spray $G^i=P(x,y)y^i$  is forward complete but not backward complete.
If $c<0$ instead, the dominant quadratic term reverses the completeness directions (the spray is backward complete but not forward complete).
\end{example}

\begin{proposition}[A nonlinear two-sided complete example]
Let $n\ge2$, let $e\in\R^n$ be a fixed unit vector, and let
$0<|\varepsilon|<1$.  
\[
\begin{array}{rl}
 P(x,y)={}&|y|+4|x|^2\langle x,y\rangle\\
 &+{\varepsilon\over2}\left\{
 \langle e,y\rangle\left(3|x|^2-{\langle x,y\rangle^2\over |y|^2}\right)\frac{\langle x,y\rangle}{|y|}
 +3\langle e,x\rangle\left(|x|^2+{\langle x,y\rangle^2\over |y|^2}\right)|y|\right\}.
\end{array}
\]
The projective spray $G^i=P(x,y)y^i$ on $\R^n$ is smooth for
$y\ne0$, satisfies the Hamel equation, is genuinely nonlinear and
nonflat, and is geodesically complete in both affine-time directions.
\end{proposition}

\noindent{\it Proof.}
Choose the homogeneous potentials
\[
 \Phi_0(y)=|y|,\qquad
 \Phi_3(y)=|y|^4+\varepsilon |y|^3\langle e,y\rangle,
 \qquad \Phi_m=0\quad(m\ne0,3).
\]
They have degrees one and four, respectively.  The representation
formula yields
\[
 P(x,y)=\Phi_0(y)+{1\over6}(x^i\partial_{y^i})^3\Phi_3(y),
\]
whose explicit expansion is the stated formula.  In particular,
$P$ satisfies the Hamel equation.  It is not linear in $y$, since
at $x=0$ one has $P(0,y)=|y|$, which is not linear in $y$.
Also $P_{x^k}(0,y)=0$, so
\[
 (P^2-P_{x^k}y^k)(0,y)=|y|^2\ne0.
\]

For an arbitrary oriented line $x=x_0+sv$, $v\ne0$, the function
$p(s)=P(x_0+sv,v)$ is a polynomial of degree three.  Its leading
coefficient, by Euler's identity for the degree-four potential, is
\[
 {1\over6}(v\cdot\partial_y)^3\Phi_3(v)
 =4\Phi_3(v)
 =4|v|^4\left(1+\varepsilon{\langle e,v\rangle\over|v|}\right)
 \ge4(1-|\varepsilon|)|v|^4>0.
\]
Thus $p(s)=A(v)s^3+O(s^2)$, where $A(v)>0$, and
\[
 2\int_0^s p(\sigma)\,d\sigma
 ={A(v)\over2}s^4+O(|s|^3)\longrightarrow+\infty
 \quad\hbox{as }s\to\pm\infty.
\]
The affine-time density
$dt/ds=C\exp(2\int_0^s p(\sigma)\,d\sigma)$ is therefore
nonintegrable at both ends of every oriented line.  Every geodesic
extends for all positive and negative affine times. \hfill {\it Q.E.D.}

\noindent
Zhongmin Shen\\
Department of Mathematical Sciences\\
Indiana University Indianapolis\\
402 N Blackford Street\\
Indianapolis, IN 46202\\
USA\\

\noindent
Email: zshen@iu.edu

\end{document}